\documentclass[a4paper]{article}

\usepackage[margin=1in]{geometry}

\usepackage{epsfig}
\usepackage{amsfonts}
\usepackage{amssymb}
\usepackage{graphics}
\usepackage{mathrsfs}

\usepackage[T1]{fontenc}
\newcommand{\changefont}[3]{
\fontfamily{#1} \fontseries{#2} \fontshape{#3} \selectfont}

\changefont{ptm}{m}{n}

\usepackage{setspace} 
 \usepackage{graphicx}    

\usepackage{epsfig}
\usepackage{amsfonts}
\usepackage{amssymb}
\usepackage{graphics}
\usepackage{mathrsfs}
\newtheorem{thm}{Theorem}[section]
\newtheorem{defn}{Definition}[section]

\newtheorem{lem}{Lemma}[section]

\newtheorem{remark}{Remark}[section]

\newcommand{\rr}{\mathbb R}
\newcommand{\zz}{\mathbb Z}

\usepackage{amsfonts}
\usepackage{amssymb}
\usepackage{graphics}
\usepackage{mathrsfs}
\usepackage{color}

\long\def\symbolfootnote[#1]#2{\begingroup%
\def\thefootnote{\fnsymbol{footnote}}\footnote[#1]{#2}\endgroup} 

\begin{document}

%

\begin{center}
\Large \textbf{Almost periodic solutions of retarded SICNNs with functional response on piecewise constant argument}
\end{center}

\vspace{-0.3cm}
\begin{center}
\normalsize \textbf{Marat Akhmet$^{a,} \symbolfootnote[1]{Corresponding Author Tel.: +90 312 210 5355,  Fax: +90 312 210 2972, E-mail: marat@metu.edu.tr}$, Mehmet Onur Fen$^b$, Mokhtar Kirane$^{c}$} \\
\vspace{0.2cm}
\textit{\textbf{\footnotesize$^a$Department of Mathematics, Middle East Technical University, 06800, Ankara, Turkey}} \\
\textit{\textbf{\footnotesize$^b$Neuroscience Institute, Georgia State University, Atlanta, Georgia 30303, USA}} \\
\textit{\textbf{\footnotesize$^c$Laboratoire de Mathématiques, Image et Applications, Pôle Sciences et Technologies, Université de La Rochelle, Avenue Michel Crépeau, 17042 La Rochelle, France}}
\vspace{0.1cm}
\end{center}

\vspace{0.3cm}

\begin{center}
\textbf{Abstract}
\end{center}

\noindent\ignorespaces

We consider a new model for shunting inhibitory cellular neural networks, retarded functional differential equations with piecewise constant argument. The existence and exponential stability of almost periodic solutions are investigated. An illustrative example is provided.

\vspace{0.2cm}
 
\noindent\ignorespaces \textbf{Keywords:} Shunting inhibitory cellular neural networks; Retarded functional differential equations; Alternate constancy of argument; Bohr almost periodic solutions; Exponential stability

\section{Introduction}

Cellular neural networks ($CNNs$) have been paid much attention in the past two decades \cite{chua}-\cite{1}. Exceptional role in psychophysics, speech, perception, robotics, adaptive pattern recognition, vision, and image processing is played by shunting inhibitory cellular neural networks ($SICNNs$), which was introduced by Bouzerdoum and Pinter \cite{bouzer1}. One of the most attractive subjects for this type of neural networks is the existence of almost periodic solutions. This problem has been investigated for models with different types of activation functions \cite{10}-\cite{12}. In the present study, we investigate a new model of $SICNNs$ by considering deviated as well as piecewise constant time arguments, and prove the existence of exponentially stable almost periodic solutions. All the results are discussed for the general type of activation functions, but they can be easily specified for applications.

Extended information about differential equations with generalized piecewise constant argument \cite{a1} can be found in the book \cite{akhmet}. As a subclass, they contain differential equations with piecewise constant argument ($EPCA$) \cite{aw1}-\cite{w1}, where the piecewise constant argument is assumed to be a multiple of the greatest integer function.

Differential equations with piecewise constant argument are very useful as models for neural networks. This was shown in the studies \cite{hul}-\cite{xue}, where the authors utilized $EPCA.$ We propose to involve a new type of systems, retarded functional differential equations with piecewise constant argument of generalized type, in the modeling. It will help to investigate a larger class of neural networks.    

In paper \cite{a1}, differential equations with piecewise constant argument of generalized type ($EPCAG$) were introduced. We not only maximally generalized the argument functions, but also proposed to reduce investigation of $EPCAG$ to integral equations. Due to that innovation, it is now possible to analyze essentially non-linear systems, that is, systems non-linear with respect to values of solutions at discrete moments of time, where the argument changes its constancy. Previously, the main and unique method for $EPCA$ was reduction to discrete equations and, hence, only equations in which values of solutions at the discrete moments appear linearly \cite{aw1}-\cite{w1} have been considered. 

The crucial novelty of the present paper is that the piecewise constant argument in the functional differential equations is of alternate (advanced-delayed) type. In the literature, biological reasons for the argument to be delayed were discussed \cite{murray,peskin}. However, the role of advanced arguments has not been analyzed properly yet. Nevertheless, the importance of anticipation for biology was mentioned by some authors. For example, in the paper \cite{bucks}, it is supposed that synchronization of biological oscillators may request anticipation of counterparts behavior. Consequently, one can assume that equations for neural networks may also need anticipation, which is usually reflected in models by advanced argument. Therefore, the systems taken into account in the present study can be useful in future analyses of $SICNNs.$ Furthermore, the idea of involving both advanced and delayed arguments in neural networks can be explained by the existence of retarded and advanced actions in a model of classical electrodynamics \cite{driver}. Moreover, mixed type deviation of the argument may depend on traveling waves emergence in $CNNs$ \cite{waves}. Understanding the structure of such traveling waves is important due to their potential applications including image processing (see, for example, \cite{chua}-\cite{waves}). More detailed analysis of deviated arguments in neural networks can be found in \cite{ay}-\cite{aay2}.
 
Shunting inhibition is a phenomenon in which the cell is ``clamped'' to its resting potential when the reversal potential of $Cl^-$ channels are close to the membrane resting potential of the cell \cite{bouzer1,Shepherd04}. It occurs through the opposition of an inward current, which would otherwise depolarize the membrane potential to threshold, by an inward flow of $Cl^-$ ions \cite{Shepherd04}. From the biological point of view, shunting inhibition has an important role in the dynamics of neurons \cite{Vida06}-\cite{Graham98}. According to the results of Vida et al. \cite{Vida06} networks with shunting inhibition are advantageous compared to the networks with hyperpolarizing inhibition such that in the former type networks oscillations are generated with smaller tonic excitatory drive, network frequencies are tuned to the $\gamma$ band, and robustness against heterogeneity in the excitatory drive is markedly improved. It was demonstrated by Mitchell and Silver \cite{Mitchell03} that shunting inhibition can modulate the gain and offset of the relationship between output firing rate and input frequency in granule cells when excitation and/or inhibition are mediated by time dependent synaptic input. Besides, Borg-Graham et al. \cite{Graham98} proposed that nonlinear shunting inhibition may act during the initial stage of visual cortical processing, setting the balance between opponent `On' and `Off' responses in different locations of the visual receptive field \cite{Graham98}. On the other hand, shunting neural networks are important for various engineering applications \cite{bouzer1},\cite{bouzer2}-\cite{Arulampalam01}. For example, in vision, shunting lateral inhibition enhances edges and contrast, mediates directional selectivity, and causes adaptation of the organization of the spatial receptive field and of the contrast sensitivity function \cite{bouzer1},\cite{bouzer2}-\cite{Pinter83}. Moreover, such networks are appropriate to be used in medical diagnosis \cite{Arulampalam01}. Therefore, the investigation of the dynamics of $SICNNs,$ which are biologically inspired networks designed upon the shunting inhibition concept \cite{bouzer1}, is important for the improvement of the techniques used in medical diagnosis, adaptive pattern recognition, image processing etc. \cite{bouzer2}-\cite{Arulampalam01} and may shed light on neuronal activities concerning shunting inhibition \cite{Vida06}-\cite{Graham98}. 
 
Exponential stability of neural networks has been widely studied in the literature (see, for example, \cite{Chunxia}-\cite{9},\cite{Liao02}-\cite{7}). According to Liao et al. \cite{Liao02}, the exponential stability has importance in neural networks when the exponentially convergence rate is used to determine the speed of neural computations. The studies \cite{Liao02,Yi99} were concerned with the exponential stability and estimation of exponential convergence rates in neural networks. In the paper \cite{Liao02}, Lyapunov-Krasovskii functionals and the linear matrix inequality (LMI) approaches were combined to investigate the problem, whereas the boundedness of the Dini derivative of the neuron input output activations was required in \cite{Yi99}. The exponential stabilization problem of memristive neural networks was considered in \cite{Wen15a} by means of the Lyapunov-Krasovskii functional and free weighting matrix techniques. Additionally, the Lyapunov-Krasovskii functional method was considered by Wen et al. \cite{Wen15c} to analyze the passivity of stochastic impulsive memristor-based piecewise linear systems, and the free weighting matrix approach was utilized in \cite{He06} to derive an LMI based delay dependent exponential stability criterion for neural networks with a time varying delay. On the other hand, exponential stability criteria were derived by Dan et al. \cite{Dan13} for an error system in order to achieve lag synchronization of coupled delayed chaotic neural networks. The concept of lag synchronization was taken into account also within the scope of the papers \cite{Wen14} and \cite{Wen15b} for memristive neural networks and for a class of switched neural networks with time-varying delays, respectively. Furthermore, the Banach fixed point theorem and the variant of a certain integral inequality with explicit estimate were used to investigate the global exponential stability of pseudo almost periodic solutions of $SICNNs$ with mixed delays in the study \cite{cherif}.

Almost periodic and in particular quasi-periodic motions are important for the theory of neural networks. According to Pasemann et al. \cite{Pasemann03}, periodic and quasi-periodic solutions have many fundamental importances in biological and artificial systems, as they are associated with central pattern generators, establishing stability properties and bifurcations (leading to the discovery of periodic solutions). Besides, the sinusoidal shape of neural output signals is, in general, associated with appropriate quasi-periodic attractors for discrete-time dynamical systems. In the book \cite{Izhikevich07}, the dynamics of the brain activity is considered as a system of many coupled oscillators with different incommensurable periods. Signals from the neurons have a phase shift of $\pi/2,$ and may be useful for various kinds of applications; for instance, controlling the gait of legged robots \cite{Kimura99}. Furthermore, an alternative discrete time model of coupled quasi-periodic and chaotic neural network oscillators were considered by Wang \cite{Wang92}.

Let us describe the model of $SICNNs$ in its most original form \cite{bouzer1}. Consider a two dimensional grid of processing cells arranged into $m$ rows and $n$ columns, and let $C_{ij},$ $i=1,2,\ldots,m,$ $j=1,2,\ldots,n,$ denote the cell at the $(i,j)$ position of the lattice. In $SICNNs,$ neighboring cells exert mutual inhibitory interactions of the shunting type. The dynamics of a cell $C_{ij}$ are described by the following nonlinear ordinary differential equation,
\begin{eqnarray} 
\begin{array}{l} \label{1}
\displaystyle \frac{dx_{ij}}{dt}=-a_{ij}x_{ij}-\sum_{C_{kl}\in N_{r}(i,j)} C_{ij}^{kl}f(x_{kl}(t))x_{ij} + L_{ij}(t),
\end{array}
\end{eqnarray}
where $x_{ij}$ is the activity of the cell $C_{ij};$ $L_{ij}(t)$ is the external input to the cell $C_{ij};$ the constant $a_{ij}>0$ represents the passive decay rate of the cell activity; $C_{ij}^{kl}\geq 0$ is the coupling strength of postsynaptic activity of the cell $C_{kl}$ transmitted to the cell $C_{ij};$ the activation function $f(x_{kl})$ is a positive continuous function representing the output or firing rate of the cell $C_{kl};$ and the $r-$neighborhood of the cell $C_{ij}$ is defined as
\begin{eqnarray*} 
N_{r}(i,j)=\{C_{kl}: \max(|k-i|,|l-j|)\leq r,  \ 1\leq k\leq m,   1\leq l \leq n \}.
\end{eqnarray*} 

It is worth noting that even if the activation function is supposed to be globally bounded and Lipschitzian, these properties are not valid for the nonlinear terms in the right hand sides of the differential equations describing the dynamics of $SICNNs,$ and this is one of the reasons why a sophisticated mathematical analysis is required for $SICNNs$ in general. Another reason is that the
connections between neurons in $SICNNs$ act locally only in $r-$neighborhoods. This causes special ways of evaluations different than those customized for earlier developed neural networks in the mathematical analyses of the models.

It is reasonable to say that the usage of deviated arguments in neural networks makes the models much closer to applications. For example, in \cite{8} the model  was  considered with variable delays,
\begin{eqnarray}\label{2}
\frac{dx_{ij}}{dt} =  -  a_{ij} x_{ij} -  \sum_{C_{kl} \in N_r(i,j)} C_{ij}^{kl}f(x_{kl}(t- \tau(t)))x_{ij} +  L_{ij}(t).
\end{eqnarray}
In the present study, we introduce and investigate more general neural networks. The model will be described in the next section.

\section{Preliminaries}

Let $\zz$ and $\rr$ denote the sets of all integers and real numbers, respectively. Throughout the paper, the norm $\left\|u\right\|=\displaystyle \max_{(i,j)} \left|u_{ij}\right|,$ where $u=\left\{u_{ij}\right\} = (u_{11},\ldots,u_{1n}, \ldots, u_{m1} \ldots,u_{mn}) \in \mathbb R^{m\times n},$ will be used.

Suppose that $\theta=\{\theta_p\}$ and $\zeta=\{\zeta_p\},$ $p \in \zz,$ are sequences of real numbers such that the first one is strictly ordered, $|\theta_p| \to \infty$ as $|p| \to \infty,$ and the second one satisfies $\theta_p \le \zeta_p \le \theta_{p+1}$ for all $p \in \zz.$ The sequence $\zeta$ is not necessarily strictly ordered. We say that a function is of $\gamma-$type, and denote it by $\gamma(t),$ if $\gamma(t) = \zeta_p$ for $\theta_p \le t < \theta_{p+1},$ $p \in \zz.$ One can affirm, for example, that $\displaystyle 2 \left[\frac{t+1}{2} \right]$ is a $\gamma-$type function with $\theta_p=2p-1$, $\zeta_p=2p.$

Fix a non-negative number $\tau \in \rr$ and let ${\cal C}^0$ be the set of all continuous functions mapping the interval $[-\tau,0]$ into $\mathbb R,$ with the uniform norm $\displaystyle \|\phi\|_0 =  \max_{t\in[-\tau,0]} \left| \phi(t) \right|.$ Moreover, we denote by ${\cal C}$ the set consisting of continuous functions mapping the interval $[-\tau,0]$ into  $\rr^{m \times n},$ with the uniform norm $\displaystyle \|\phi\|_0 =  \max_{t\in[-\tau,0]}  \left\| \phi(t) \right\|.$

In the present study, we propose to investigate retarded $SICNNs$ with functional response on piecewise constant argument of the following form,  
\begin{eqnarray}\label{3}
\frac{dx_{ij}}{dt} =  - a_{ij} x_{ij} - \sum_{C_{kl} \in N_r(i,j)} C_{ij}^{kl}f(x_{klt},x_{kl\gamma(t)})x_{ij} + L_{ij}(t),
\end{eqnarray} 
where $f:{\cal C}^0 \times {\cal C}^0 \to \mathbb R$ is a continuous functional. 

In network (\ref{3}), the terms $x_{klt}$ and $x_{kl\gamma(t)}$ must be understood in the way used for functional differential equations \cite{burton}-\cite{kuang}. That is, $x_{klt}(s) = x_{kl}(t+s)$ and $x_{kl\gamma(t)}(s) = x_{kl}(\gamma(t) + s)$ for $s \in [-\tau,0].$  Let us clarify that the argument function $\gamma(t)$ is of the alternate type. Fix an integer $p$ and consider the  function  on  the interval $[\theta_{p},\theta_{p+1}).$ Then, the function $\gamma(t)$ is equal to  $\zeta_{p}.$ If the argument $t$ satisfies $\theta_p \leq t < \zeta_p,$ then $\gamma(t)> t$ and it is of advanced type. Similarly, if $\zeta_p < t < \theta_{p+1},$ then $\gamma(t)< t$ and, hence, it is of delayed type. Consequently, it is worth noting that the $SICNN$ (\ref{3}) is with {\it alternate constancy} of argument. It is known that $\gamma(t)$ is the most general among piecewise constant argument functions \cite{akhmet}. Our model is much more general than the equations investigated in \cite{s3}-\cite{wyz1}, where the delay is constant $\tau = 1$ and it is equal to the step of the greatest integer function $[t].$ Differential equations with functional response on the piecewise constant argument were first introduced in the paper \cite{a9}. In the present study, we apply the theory to the analysis of neural networks. All previous authors were at most busy with terms of the form $x({\gamma(t)}).$ Thus, one can say that retarded functional differential equations with piecewise constant argument in the most general form is investigated in this paper.

Since the model (\ref{3}) is a new one, we have to investigate not only the existence of almost periodic solutions and their stability, but also common problems of  the existence and uniqueness of solutions, their continuation to infinity and boundedness.

One can easily see that system (\ref{2}) is a particular case of (\ref{3}). Additionally, results of the present paper are true or can be easily adapted to the following systems, 
\begin{eqnarray}\label{4}
\frac{dx_{ij}}{dt} =  -  a_{ij} x_{ij} -  \sum_{C_{kl} \in N_r(i,j)} C_{ij}^{kl}f(x_{kl}(\gamma(t)))x_{ij} +  L_{ij}(t),
\end{eqnarray} 
that is,   differential  equations with  piecewise constant  argument, $EPCAG,$
\begin{eqnarray}\label{5}
\frac{dx_{ij}}{dt} =  -  a_{ij} x_{ij} -  \sum_{C_{kl} \in N_r(i,j)} C_{ij}^{kl}f(x_{kl}(t - \tau(t)),x_{kl}(\gamma(t)- \tau(t)))x_{ij} +  L_{ij}(t),
\end{eqnarray}
differential  equations with  variable delay and piecewise  constant  argument, 
\begin{eqnarray}\label{6}
\frac{dx_{ij}}{dt} =  -  a_{ij} x_{ij} -  \sum_{C_{kl} \in N_r(i,j)} C_{ij}^{kl}f(x_{kl}(t - \tau(t)))x_{ij} -   \sum_{C_{kl} \in N_r(i,j)} D_{ij}^{kl}g(x_{kl}(\gamma(t)))x_{ij} + L_{ij}(t).
\end{eqnarray}  
In other words, what we have suggested are sufficiently general models, which can be easily specified for concrete applications.

Let us introduce the initial condition for $SICNN$ (\ref{3}). Fix a number $\sigma \in \rr$ and functions $\phi=\left\{\phi_{ij}\right\},$ $\psi=\left\{\psi_{ij} \right\} \in {\cal C},$ $i=1,2,\ldots,m,$ $j=1,2,\ldots,n.$ In the case  $\gamma(\sigma) < \sigma,$ we say that a solution $x(t)= \left\{x_{ij}(t)\right\}$ of (\ref{3}) satisfies the initial condition and write $x(t) = x(t,\sigma,\phi,\psi),$ $t \ge \sigma,$ if $x_{\sigma}(s) = \phi(s),$ $x_{\gamma(\sigma)}(s) = \psi(s)$ for $s \in [-\tau,0].$ In what follows, we assume that if the set $[\gamma(\sigma)-\tau,\gamma(\sigma)] \cup [\sigma -\tau,\sigma]$ is connected, then the equation $\phi(s) = \psi(s+\sigma - \gamma(\sigma))$ is true for  all $s \in [-\tau, \gamma(\sigma) - \sigma].$ 
If  $\gamma(\sigma) \ge \sigma,$ then we look for a solution $x(t)= x(t,\sigma,\phi),$ $t \ge \sigma,$ such  that  $x_{\sigma}(s)  = \phi(s), s \in [-\tau,0].$
 Thus, if $\theta_p \le \sigma < \theta_{p+1}$ for some $p \in \zz,$ then there are two cases of the initial  condition:
\begin{itemize}
\item[\bf (IC$_1$)] $x_{\sigma}(s) = \phi(s),$ $\phi \in {\cal C},$ $s \in [-\tau,0]$ if $\theta_p \le \sigma  \le \zeta_p < \theta_{p+1};$
\item[\bf (IC$_2$)] $x_{\sigma}(s) = \phi(s),$ $x_{\gamma(\sigma)}(s) = \psi(s),$ $\phi,\psi \in {\cal C},$ $s \in [-\tau,0],$ if  $\theta_p \le \zeta_p < \sigma < \theta_{p+1}.$	
\end{itemize}
Considering $SICNN$ (\ref{3}) with these conditions, we shall say about the initial value problem $(IVP)$ for (\ref{3}). To be short, we shall say only about $IVP$ in the form  $x(t,\sigma,\phi,\psi),$ specifying $x(t,\sigma,\phi)$ for $(IC_1),$ if needed. Thus, we can provide the following definitions now. 

\begin{defn}\label{d2}
A function $x(t)= \left\{x_{ij}(t)\right\},$ $i=1,2,\ldots,m,$ $j=1,2,\ldots,n,$ is a solution of (\ref{3}) with $(IC_1)$ or $(IC_2)$ on an interval $[\sigma, \sigma + a)$ if:
\begin{enumerate}
\item  [(i)] it satisfies the initial condition;
\item[(ii)] $x(t)$ is continuous on $[\sigma,  \sigma + a);$
\item[(iii)] the derivative $x'(t)$ exists for  
$t\geq \sigma$ with the possible exception of the points
$\theta_p$,  where one-sided derivatives exist;  
\item[(iv)] equation $(\ref{3})$ is satisfied by $x(t)$ for all $t > \sigma$ except possibly at the points of $\theta,$ and it holds for the right derivative of $x(t)$ at the points $\theta_p.$
\end{enumerate}
\end{defn}

\begin{defn}\label{defn1}
A function $x(t)= \left\{x_{ij}(t)\right\},$ $i=1,2,\ldots,m,$ $j=1,2,\ldots,n,$ is a solution of (\ref{3}) on $\rr$ if:

\begin{enumerate}

\item[(i)] $x(t)$ is continuous;

\item[(ii)] the derivative $x'(t)$ exists for all 
$t \in \rr$ with the possible exception of the points
$\theta_p,$ $p \in \zz,$  where one-sided derivatives exist;  
\item[(iv)] equation $(\ref{3})$ is satisfied by $x(t)$ for all $t \in \rr$ except at the points of $\theta,$ and it holds for the right derivative of $x(t)$ at the points $\theta_p,$ $p \in \zz.$
\end{enumerate}
\end{defn}
 
The existence and uniqueness of solutions of (\ref{3}) will be investigated in the next section.
  
\section{Existence and uniqueness}

Throughout the paper we suppose in $SICNN$ (\ref{3}) that $\displaystyle \gamma_0 = \min_{(i,j)} a_{ij} > 0$ and $ C_{ij}^{kl}$ are non-negative numbers. 
 
The following assumptions are required.
\begin{itemize}
\item[\bf (C1)] The functional $f$ satisfies the Lipschitz condition
\[
\|f(\phi_1,\psi_1) - f(\phi_2,\psi_2)\|  \le L( \|\phi_1 - \phi_2\|_0 +  \|\psi_1 - \psi_2\|_0 ),
\]
 for some positive constant $L,$ where $(\phi_1,\psi_1)$ and $(\phi_2,\psi_2)$ are from ${\cal C}^0 \times {\cal C}^0;$
\item [\bf (C2)] There exists a positive number $M$ such that $\displaystyle \sup_{(\phi,\psi) \in {\cal C}^0\times {\cal C}^0} \left|f (\phi,\psi)\right| \le M;$ 
\item [\bf (C3)] There exists a positive number $\bar \theta$ such that $\theta_{p+1} - \theta_{p} \leq \bar \theta$ for all $p \in \mathbb Z;$ 
\item [\bf (C4)] $|L_{ij}(t)| \le L_{ij}$ for all $i,j$ and $t \in \rr,$ where $L_{ij}$ are non-negative real constants.   
\end{itemize}

In the remaining parts of the paper, the notations
$\mu = \displaystyle \max_{(i,j)}\sum_{C_{kl} \in N_r(i,j)} C_{ij}^{kl},$ $\displaystyle \bar c  =  \max_{(i,j)}\frac{\sum_{C_{kl} \in N_r(i,j)} C_{ij}^{kl}}{a_{ij}},$ $\displaystyle \bar d  =  \max_{(i,j)}\frac{\sum_{C_{kl} \in N_r(i,j)} C_{ij}^{kl}}{2a_{ij}-\gamma_0},$ $\displaystyle \bar L= \max_{(i,j)} L_{ij}$ and $\displaystyle \bar l = \max_{(i,j)}\frac{L_{ij}}{a_{ij}}$ will be used. We assume that $\mu \bar{\theta}M<1$ and $M \bar c<1.$ 

Let us denote ${\cal C}_{H_0} = \{\phi \in {\cal C} : \|\phi\|_0\leq H_0\},$ where  $H_0$ is a positive number.

\begin{lem} \label{tip-tap} Suppose that the conditions $(C1)-(C4)$ hold and fix an integer $p.$ If $H_0$ is a positive number such that $\displaystyle \mu \bar{\theta} \Big[ M + \frac{2L(H_0+\bar{\theta} \bar{L})}{1-\mu \bar{\theta}M} \Big]<1,$ then for every $(\sigma,\phi,\psi) \in [\theta_p, \theta_{p+1}] \times {\cal C}_{H_0} \times {\cal C}_{H_0}$ there exists a unique solution $x(t)=x(t,\sigma, \phi,\psi)$ of (\ref{3}) on $[\sigma, \theta_{p+1}].$ 
\end{lem}

\noindent {\bf Proof.} We assume without loss of generality that $\theta_p\leq \sigma \le \zeta_p < \theta_{p+1}.$ That is, we  consider $(IC_1)$ and the solution $x(t,\sigma, \phi).$

Fix an arbitrary function $\phi \in {\cal C}_{H_0}.$ Let us denote by $\Lambda$ the set of continuous functions $u(t)=\left\{u_{ij}(t)\right\},$ $i=1,2,\ldots,m,$ $j=1,2,\ldots,n,$ defined on $[\sigma-\tau,\theta_{p+1}]$ such that $u_{\sigma}(t)=\phi(t),$ $t\in [-\tau,0],$ and $\left\|u\right\|_1 \le K_0,$ where $\displaystyle \left\|u\right\|_1 = \max_{t\in [\sigma,\theta_{p+1}]} \left\|u(t)\right\|$ and $K_0=\displaystyle \frac{H_0+\bar{\theta}\bar L}{1-\mu \bar{\theta}M}.$

Define on $\Lambda$ an operator $\Phi$ such that 
\begin{eqnarray*}
(\Phi u (t))_{ij}= \left\{\begin{array}{ll} \phi_{ij}(t-\sigma), ~t \in [\sigma-\tau,\sigma],\\
 e^{-a_{ij}(t-\sigma)}\phi_{ij}(0) - \displaystyle \int_{\sigma}^{t}e^{-a_{ij}(t-s)}\Big[\sum_{C_{kl} \in N_r(i,j)} C_{ij}^{kl}\\
\times f(u_{kls}, u_{kl\gamma(s)} )u_{ij}(s) -L_{ij}(s) \Big] ds, ~t \in [\sigma, \theta_{p+1}].
\end{array}\right.
\end{eqnarray*}

One can confirm that
$
\left| (\Phi u (t))_{ij}  \right| \le H_0 + \Big( MK_0 \displaystyle \sum_{C_{kl} \in N_r(i,j)} C_{ij}^{kl} + \bar{L} \Big) \bar{\theta},
$
$t\in [\sigma,\theta_{p+1}].$
Accordingly, the inequality
$
\left\|\Phi u\right\|_1 \le H_0 + (\mu M K_0 + \bar{L})\bar{\theta} = K_0
$
is valid. Therefore, $\Phi(\Lambda) \subseteq \Lambda.$

On the other hand, if $u(t)=\left\{u_{ij}(t)\right\}$ and $v(t)=\left\{v_{ij}(t)\right\}$ belong to $\Lambda,$ then we have for $t\in [\sigma,\theta_{p+1}]$ that 
\begin{eqnarray*}
&& \left| (\Phi u(t))_{ij} - (\Phi v(t))_{ij} \right| \le  \displaystyle \int_{\sigma}^{t}e^{-a_{ij}(t-s)} \sum_{C_{kl} \in N_r(i,j)} C_{ij}^{kl} \left|f(u_{kls},u_{kl\gamma(s)})\right| \left|u_{ij}(s)-v_{ij}(s)\right| ds \\
&& + \displaystyle \int_{\sigma}^{t}e^{-a_{ij}(t-s)} \sum_{C_{kl} \in N_r(i,j)} C_{ij}^{kl}  \left| f(u_{kls},u_{kl\gamma(s)}) - f(v_{kls},v_{kl\gamma(s)})\right|  \left|v_{ij}(s)\right|  ds \\
&& \le \bar{\theta} (M+2K_0L) \left\|u-v\right\|_1  \sum_{C_{kl} \in N_r(i,j)} C_{ij}^{kl}.
\end{eqnarray*}
Hence, the inequality $\left\|\Phi u - \Phi  v \right\|_1 \le \mu \bar{\theta} (M+2K_0L) \left\|u-v\right\|_1$ holds. Because $\mu \bar{\theta} (M+2K_0L) = \displaystyle \mu \bar{\theta} \Big[ M + \frac{2L(H_0+\bar{\theta} \bar{L})}{1-\mu \bar{\theta}M} \Big] <1,$  the operator $\Phi$ is a contraction.
Consequently, there exists a unique solution of (\ref{3}) on $[\sigma, \theta_{p+1}].$ $\square$
 
The  next assertion can  be proved  exactly  in the way  that is used to  verify Lemma $2.2$  from \cite{akhmet}, if we use Lemma \ref{tip-tap}. 
\begin{lem}\label{lemi2} 
Suppose that the conditions $(C1)-(C4)$ hold and fix an integer $p.$ If $H_0$ is a positive number such that $\displaystyle \mu \bar{\theta} \Big[ M + \frac{2L(H_0+\bar{\theta} \bar{L})}{1-\mu \bar{\theta}M} \Big]<1,$ then for every $(\sigma,\phi,\psi) \in [\theta_p, \theta_{p+1}] \times {\cal C}_{H_0} \times {\cal C}_{H_0}$ there exists a unique solution $x(t)=x(t,\sigma, \phi,\psi),$ $t \ge \sigma,$ of (\ref{3}), and it   satisfies  the integral equation
\begin{eqnarray}\label{pshik2}
x_{ij}(t)= {\rm e}^{-a_{ij}(t-\sigma)}\phi_{ij}(0) - \int_{\sigma}^{t}{\rm e}^{-a_{ij}(t-s)}\Big[\sum_{C_{kl} \in N_r(i,j)} C_{ij}^{kl}f(x_{kls},x_{kl\gamma(s)})x_{ij}(s)  -L_{ij}(s)\Big]ds.
\end{eqnarray}
\end{lem}

\section{Bounded solutions}

In this section, we will investigate the existence of a unique bounded solution of $SICNN$ (\ref{3}). Moreover, the exponential stability of the bounded solution will be considered. An auxiliary result is presented in the following lemma.

\begin{lem}\label{lemi3} 
Assume that the conditions $(C1)-(C4)$ are fulfilled. If $H_0$ is a positive number such that $\displaystyle \mu \bar{\theta} \Big[ M + \frac{2L(H_0+\bar{\theta} \bar{L})}{1-\mu \bar{\theta}M} \Big]<1,$ then a function $x(t)=\left\{x_{ij}(t)\right\},$ $i=1,2,\ldots,m,$ $j=1,2,\ldots,n,$ satisfying $\displaystyle\sup_{t\in\mathbb R} \left\|x(t)\right\|\le H_0$ is a solution of (\ref{3}) if and only if it satisfies the following integral equation
\begin{eqnarray}\label{pshik3}
&& x_{ij}(t)=  - \int_{-\infty}^{t}{\rm e}^{-a_{ij}(t-s)}\Big[\sum_{C_{kl} \in N_r(i,j)} C_{ij}^{kl}f(x_{kls},x_{kl\gamma(s)})x_{ij}(s)  -L_{ij}(s)\Big]ds.
\end{eqnarray} 
\end{lem}

\noindent {\bf Proof.} We consider only sufficiency. The necessity can be proved by using (\ref{pshik2}) in a very similar way to the ordinary differential equations case.  
One can obtain that  
\begin{eqnarray*}   
&& \displaystyle  \Big| \int_{-\infty}^{t}{\rm e}^{-a_{ij}(t-s)}\Big[\sum_{C_{kl} \in N_r(i,j)} C_{ij}^{kl}f(x_{kls},x_{kl\gamma(s)})x_{ij}(s) -L_{ij}(s)\Big]ds   \Big|\\
&& \le  \frac{1}{a_{ij}}  \Big(\displaystyle M  H_0 \sum_{C_{kl} \in N_r(i,j)} C_{ij}^{kl} + L_{ij} \Big). 
\end{eqnarray*}  
Therefore, the integral in (\ref{pshik3}) is convergent. 
Differentiate  (\ref{pshik3}) to verify that it is a solution of (\ref{3}). 
$\square$

The following conditions are needed.
\begin{itemize}
\item [\bf (C5)] $\displaystyle \mu \bar{\theta}\Big[  M + \frac{2L(H+\bar{\theta} \bar{L})}{1-\mu \bar{\theta}M} \Big]<1,$ where $H=\displaystyle \frac{\bar l}{1- M \bar c};$ 
\item [\bf (C6)] $(M+2LH)\bar c<1;$   
\item [\bf (C7)] $2\bar{d} \left[M+LHe^{\gamma_0\tau/2} \left(1+e^{\gamma_0 \bar{\theta}/2}\right)\right]<1.$

\end{itemize}

The main result concerning the existence and exponential stability of bounded solutions of (\ref{3}) is mentioned in the next theorem.

\begin{thm} \label{thm2} Suppose that the conditions $(C1)-(C6)$ hold. Then, (\ref{3}) admits a unique bounded on $\rr$ solution, which satisfies (\ref{pshik3}). If, additionally, the condition $(C7)$ is valid, then the solution is exponentially stable with exponential convergence rate $\gamma_0/2.$ 
\end{thm}

\noindent {\bf Proof.} Let $C_0(\rr)$ be the set of uniformly continuous functions defined on $\rr$ such that if $u(t) \in C_0(\rr),$ then $\|u\|_{\infty} \le H,$ where $\|u\|_{\infty} = \displaystyle \sup_{t\in\rr}\|u(t)\|.$ Define on $C_0(\rr)$ the operator $\Pi$ as
\begin{eqnarray} \label{theorem_operator_defn}
(\Pi u(t))_{ij}  \equiv  - \int_{-\infty}^{t}{\rm e}^{-a_{ij}(t-s)}\Big[\sum_{C_{kl} \in N_r(i,j)} C_{ij}^{kl} f(u_{kls}, u_{kl\gamma(s)})u_{ij}(s)  - L_{ij}(s)\Big]ds.
\end{eqnarray}

If $u(t)=\left\{u_{ij}(t)\right\},$ $i=1,2,\ldots,m,$ $j=1,2,\ldots,n,$ belongs to $C_0(\rr),$ then we have that
\begin{eqnarray*}
\left|(\Pi u(t))_{ij} \right| \le \displaystyle \int_{-\infty}^t e^{-a_{ij}(t-s)} \Big( \sum_{C_{kl} \in N_r(i,j)} C_{ij}^{kl} M H + L_{ij} \Big)ds =\frac{1}{a_{ij}} \Big( \sum_{C_{kl} \in N_r(i,j)} C_{ij}^{kl} M H + L_{ij} \Big). 
\end{eqnarray*}
Utilizing the last inequality one can show that
$
 \left\|(\Pi u) \right\|_{\infty} \le  \bar c MH+\bar l=H. 
$
Therefore, $\Pi u(t)\in C_0(\mathbb R).$

Let us verify that this operator is contractive. Indeed, if $u(t)=\left\{u_{ij}(t)\right\}$ and $v(t)=\left\{v_{ij}(t)\right\}$ belong to $C_0(\rr),$ then 
\begin{eqnarray*}
&& \left| (\Pi u(t))_{ij} -(\Pi v(t))_{ij} \right| \le \displaystyle \int_{-\infty}^t e^{-a_{ij}(t-s)} \sum_{C_{kl} \in N_r(i,j)} C_{ij}^{kl} \left| f(u_{kls},u_{kl\gamma(s)})  \right| \left|u_{ij}(s)-v_{ij}(s)\right| ds \\
&& + \displaystyle \int_{-\infty}^t e^{-a_{ij}(t-s)} \sum_{C_{kl} \in N_r(i,j)} C_{ij}^{kl} \left| f(u_{kls},u_{kl\gamma(s)}) - f(v_{kls},v_{kl\gamma(s)}) \right|  \left|v_{ij}(s)\right| ds \\
&& \le \displaystyle \int_{-\infty}^t e^{-a_{ij}(t-s)} \sum_{C_{kl} \in N_r(i,j)} C_{ij}^{kl} M \left|u_{ij}(s)-v_{ij}(s)\right| ds \\
&& + \displaystyle \int_{-\infty}^t e^{-a_{ij}(t-s)} \sum_{C_{kl} \in N_r(i,j)} C_{ij}^{kl}  HL \left( \left\|u_{kls}-v_{kls}\right\|_0 + \left\|u_{kl\gamma(s)}-v_{kl\gamma(s)}\right\|_0 \right) ds\\
&& \le (M+2LH) \frac{\sum_{C_{kl} \in N_r(i,j)} C_{ij}^{kl}}{a_{ij}} \left\|u-v\right\|_{\infty}.
\end{eqnarray*}
Hence, the inequality $ \left\|\Pi u-\Pi v\right\|_{\infty}\le (M+2LH)\bar c \left\|u-v\right\|_{\infty}$ is valid.
In accordance with condition $(C6),$ the operator $\Pi$ is contractive. Consequently, $SICNN$ (\ref{3}) admits a unique solution $\widetilde v(t)=\left\{\widetilde v_{ij}(t)\right\}$ that belongs to $C_0(\rr).$

We will continue with the investigation of the exponential stability. Fix an arbitrary number $\epsilon>0$ and let $\delta$ be a sufficiently small positive number such that  $\alpha_1<1,$ $\alpha_2<1,$ $\alpha_3<1$ and ${\cal K}(\delta) < \epsilon,$ where 
$\displaystyle {\cal K}(\delta) = \frac{\delta}{1 - 2\bar  d [M + LHe^{\gamma_0\tau/2}(1+e^{\gamma_0 \bar{\theta}/2})]},$ $\alpha_1=\displaystyle \mu \bar{\theta}\Big[  M + \frac{2L(H+\delta+\bar{\theta} \bar{L})}{1-\mu \bar{\theta}M} \Big],$
$\alpha_2=(M+2LH)\bar c + 4L \bar d {\cal K}(\delta)$ and
$\alpha_3=\mu \bar{\theta}(M+2LH) + 2 \mu \bar{\theta} L{\cal K}(\delta).$

Suppose that $\widetilde v_{\sigma}(s)=\eta(s),$ $s\in [-\tau,0].$ Let $u(t)=\left\{u_{ij}(t)\right\}$ be a solution of the network (\ref{3}) with $u_{\sigma}(s)=\phi(s),$ $s\in [-\tau,0],$ where the function $\phi$ satisfies the inequality $\|\phi - \eta\|_0 < \delta.$  Without loss of generality we assume that $\gamma(\sigma) \ge \sigma.$ Using Lemma \ref{lemi2} one can verify for $t\ge \sigma$ that  
\begin{eqnarray*} 
&& u_{ij}(t)  - \widetilde v_{ij}(t) = {\rm e}^{-a_{ij}(t-\sigma)} \left( \phi_{ij}(0)-\eta_{ij}(0) \right)\\
&& - \int_{\sigma}^{t}{\rm e}^{-a_{ij}(t-s)} \sum_{C_{kl} \in N_r(i,j)} C_{ij}^{kl} \Big[ f(u_{kls},u_{kl\gamma(s)})u_{ij}(s)- f(\widetilde v_{kls},\widetilde v_{kl\gamma(s)})\widetilde v_{ij}(s) \Big] ds. 
\end{eqnarray*}
Denote by  $w(t)=\left\{w_{ij}(t)\right\},$   the difference $u(t) - \widetilde v(t).$ Then, $w(t)$ satisfies the relation 
 \begin{eqnarray} \label{Shunting_almost_periodic_stability_proof}
 && w_{ij}(t)  =  {\rm e}^{-a_{ij}(t-\sigma)} \left(\phi_{ij}(0)-\eta_{ij}(0)\right)  - \int_{\sigma}^{t}{\rm e}^{-a_{ij}(t-s)}\sum_{C_{kl} \in N_r(i,j)} C_{ij}^{kl}\Big[ f(\widetilde v_{kls} +\nonumber\\
 && w_{kls}, \widetilde v_{kl\gamma(s)}+w_{kl\gamma(s)} )(\widetilde  v_{ij}(s) +w_{ij}(s)) - f(\widetilde  v_{kls}, \widetilde  v_{kl\gamma(s)}) \widetilde  v_{ij}(s) \Big]ds . 
\end{eqnarray}

We will consider equation (\ref{Shunting_almost_periodic_stability_proof}) for $\sigma =0.$  Let $\Psi_{\delta}$ be the set of all continuous functions $w(t)=\left\{w_{ij}(t)\right\}$ which are defined on $[-\tau, \infty)$ such that: 
\begin{itemize}
\item [(i)]  $w(t) = \phi(t) - \eta(t),$ $t\in[-\tau,0];$
\item [(ii)] $w(t)$ is  uniformly continuous  on  $[0,+ \infty);$
\item [(iii)]$|| w(t)|| \leq {\cal K}(\delta) e^{-\gamma_0 t/2}$ for $t\geq 0.$ 
\end{itemize}

Define on $\Psi_{\delta}$ an operator $\tilde \Pi$ such that 
\begin{eqnarray*}
&& (\tilde \Pi w (t))_{ij}= \left\{\begin{array}{ll} \phi_{ij}(t) - \eta_{ij}(t), ~t \in [-\tau,0],\\
 e^{-a_{ij}t}(\phi_{ij}(0)-\eta_{ij}(0)) - \displaystyle \int_{0}^{t}e^{-a_{ij}(t-s)}\sum_{C_{kl} \in N_r(i,j)} C_{ij}^{kl}\Big[f(\widetilde v_{kls} +\\
 w_{kls},\widetilde v_{kl\gamma(s)}+w_{kl\gamma(s)} )(\widetilde v_{ij}(s) + w_{ij}(s)) - f(\widetilde v_{kls},\widetilde v_{kl\gamma(s)})\widetilde v_{ij}(s)\Big]ds, ~t > 0.
\end{array}\right.
\end{eqnarray*}
We shall show that $\tilde \Pi : \Psi_{\delta} \rightarrow \Psi_{\delta}.$   
Indeed, it is true for $t\geq 0$ that 
\begin{eqnarray*}
&& |(\tilde \Pi w (t))_{ij}| \leq e^{-a_{ij}t}\delta +  \int_{0}^{t}  e^{-a_{ij}(t-s)}\sum_{C_{kl} \in N_r(i,j)} C_{ij}^{kl} \\
&& \times \left|  f(\widetilde v_{kls} +  w_{kls},\widetilde v_{kl\gamma(s)}+w_{kl\gamma(s)} ) - f(\widetilde v_{kls},\widetilde v_{kl\gamma(s)}) \right|  \left|\widetilde v_{ij}(s)\right| ds \\
&& + \int_{0}^{t}  e^{-a_{ij}(t-s)}\sum_{C_{kl} \in N_r(i,j)} C_{ij}^{kl} \left|f(\widetilde v_{kls} + w_{kls},\widetilde v_{kl\gamma(s)}+w_{kl\gamma(s)} )\right| \left|w_{ij}(s)\right| ds \\
&& \le e^{-a_{ij}t}\delta + \int_{0}^{t}  e^{-a_{ij}(t-s)}\sum_{C_{kl} \in N_r(i,j)} C_{ij}^{kl} H L \left( \left\| w_{kls} \right\|_0 + \left\| w_{kl\gamma(s)}  \right\|_0 \right) ds \\
&& + \int_{0}^{t}  e^{-a_{ij}(t-s)}\sum_{C_{kl} \in N_r(i,j)} C_{ij}^{kl} M {\cal K}(\delta) e^{-\gamma_0 s/ 2}ds \\
&& \le e^{-a_{ij}t}\delta + \int_{0}^{t}  e^{-a_{ij}(t-s)}\sum_{C_{kl} \in N_r(i,j)} C_{ij}^{kl} H L {\cal K}(\delta) \left( e^{\gamma_0 \tau /2} +e^{\gamma_0(\bar{\theta}+\tau)/2} \right) e^{-\gamma_0 s/ 2} ds \\
&& + \int_{0}^{t}  e^{-a_{ij}(t-s)}\sum_{C_{kl} \in N_r(i,j)} C_{ij}^{kl} M {\cal K}(\delta) e^{-\gamma_0 s/ 2}ds\\
&& =  e^{-a_{ij}t}\delta + \left(\frac{2\sum_{C_{kl} \in N_r(i,j)} C_{ij}^{kl}}{2a_{ij}-\gamma_0}\right) {\cal K}(\delta) [M + LHe^{\gamma_0\tau/2}(1+e^{\gamma_0 \bar{\theta}/2})] e^{-\gamma_0 t/ 2}.
\end{eqnarray*}
Thus, the inequality
\begin{eqnarray*}
\left\|\tilde \Pi w (t) \right\| \leq e^{-\gamma_0 t} \delta + 2 \bar d {\cal K}(\delta) [M + LHe^{\gamma_0\tau/2}(1+e^{\gamma_0 \bar{\theta}/2})] e^{-\gamma_0 t/ 2} \le {\cal K}(\delta)e^{-\gamma_0 t/ 2}
\end{eqnarray*}
is valid for $t\ge 0.$ 

Now, let $w^1(t)=\left\{w^1_{ij}(t)\right\},$ $w^2(t)=\left\{w^2_{ij}(t)\right\}$ be elements of $\Psi_{\delta}.$ One can confirm for $t\ge 0$ that
\begin{eqnarray*}
&& \left\| (\tilde \Pi w^1 (t))_{ij} - (\tilde \Pi w^2 (t))_{ij} \right\| \le \displaystyle \int_0^t e^{-a_{ij}(t-s)} \sum_{C_{kl} \in N_r(i,j)} C_{ij}^{kl} \left(\left| \widetilde{v}_{ij}(s) \right| + \left| w^2_{ij}(s) \right| \right)\\
&& \times \left|  f(\widetilde v_{kls} + w^1_{kls}, \widetilde v_{kl\gamma(s)}+w^1_{kl\gamma(s)}) - f(\widetilde v_{kls} +  w^2_{kls}, \widetilde v_{kl\gamma(s)} + w^2_{kl\gamma(s)}) \right| ds \\
&& + \displaystyle \int_0^t e^{-a_{ij}(t-s)} \sum_{C_{kl} \in N_r(i,j)} C_{ij}^{kl} \left| f(\widetilde v_{kls} + w^1_{kls}, \widetilde v_{kl\gamma(s)}+w^1_{kl\gamma(s)}) \right| \left| w^1_{ij}(s)-w^2_{ij}(s) \right| ds \\
&& \le \displaystyle \int_0^t e^{-a_{ij}(t-s)} \sum_{C_{kl} \in N_r(i,j)} C_{ij}^{kl}  L\Big( H +  {\cal K}(\delta)e^{-\gamma_0 s /2} \Big)  \Big(  \left\| w^1_{kls}-w^2_{kls}  \right\|_0  + \left\| w^1_{kl\gamma(s)}-w^2_{kl\gamma(s)}  \right\|_0 \Big) ds \\
&& + \displaystyle \int_0^t e^{-a_{ij}(t-s)} \sum_{C_{kl} \in N_r(i,j)} C_{ij}^{kl} M \left| w^1_{ij}(s)-w^2_{ij}(s) \right| ds \\
&& \le (M+2LH) \displaystyle \sup_{t \ge 0} \left\| w^1(t)-w^2(t) \right\| \frac{\sum_{C_{kl} \in N_r(i,j)} C_{ij}^{kl}}{a_{ij}} \left(1-e^{-a_{ij}t}\right) \\
&& + 4L{\cal K}(\delta) \displaystyle \sup_{t \ge 0} \left\| w^1(t)-w^2(t) \right\| \frac{\sum_{C_{kl} \in N_r(i,j)} C_{ij}^{kl}}{2a_{ij}-\gamma_0} \left( e^{-\gamma_0t/2} -e^{-a_{ij}t} \right).
\end{eqnarray*} 
Therefore, we have that
$\displaystyle \sup_{t\geq 0}||\tilde \Pi w^1 (t) - \tilde \Pi w^2 (t) || \leq \alpha_2 \sup_{t\geq 0}||w^1 (t) - w^2 (t) ||.$
Since $\alpha_2<1,$ one can conclude by using a contraction mapping argument that there exists a unique fixed point $\widetilde w(t)=\left\{\widetilde w_{ij}(t)\right\}$ of the operator $\tilde \Pi :\ \Psi_{\delta} \rightarrow \Psi_{\delta},$ which is a solution of (\ref{Shunting_almost_periodic_stability_proof}).

To complete the proof, we need to show that there does not exist a solution of (\ref{Shunting_almost_periodic_stability_proof}) with $\sigma=0$ different from $\widetilde w(t).$ Suppose that $\theta_p \le 0 < \theta_{p+1}$ for some $p \in \zz.$ Assume that there exists a solution $\overline w(t)=\left\{\overline w_{ij}(t)\right\}$ of (\ref{Shunting_almost_periodic_stability_proof}) different from $\widetilde w(t).$ Denote by $z(t)=\left\{z_{ij}(t)\right\}$ the difference $\overline w(t)-\widetilde w(t),$ and let $\displaystyle \max_{t\in [0, \theta_{p+1}]} ||z(t)|| = \bar m.$ 
It can be verified for $t\in[0,\theta_{p+1}]$ that 
\begin{eqnarray*}
&& |z_{ij}(t)| \le   \displaystyle \int_0^t e^{-a_{ij}(t-s)} \sum_{C_{kl} \in N_r(i,j)} C_{ij}^{kl} \left| \widetilde v_{ij}(s) + \widetilde w_{ij}(s) \right| \\
&& \times \left| f(\widetilde{v}_{kls}+\overline{w}_{kls}, \widetilde{v}_{kl\gamma(s)}+\overline{w}_{kl\gamma(s)}) - f(\widetilde{v}_{kls}+\widetilde{w}_{kls}, \widetilde{v}_{kl\gamma(s)}+\widetilde{w}_{kl\gamma(s)})  \right| ds \\
&& + \displaystyle \int_0^t e^{-a_{ij}(t-s)} \sum_{C_{kl} \in N_r(i,j)} C_{ij}^{kl}  \left|f(\widetilde{v}_{kls}+\overline{w}_{kls}, \widetilde{v}_{kl\gamma(s)}+\overline{w}_{kl\gamma(s)})\right|  \left| z_{ij}(s)  \right| ds \\
&& \le \displaystyle \int_0^t e^{-a_{ij}(t-s)} \sum_{C_{kl} \in N_r(i,j)} C_{ij}^{kl} L \left( H + {\cal K}(\delta) \right) \left(  \left\|  z_{kls} \right\|_0 + \left\|  z_{kl\gamma(s)} \right\|_0  \right) ds \\
&& + \displaystyle \int_0^t e^{-a_{ij}(t-s)} \sum_{C_{kl} \in N_r(i,j)} C_{ij}^{kl}  M  \left| z_{ij}(s)  \right| ds \\
&& \le \bar{\theta} \bar{m}\left[M+2L \left( H + {\cal K}(\delta) \right)\right] \sum_{C_{kl} \in N_r(i,j)} C_{ij}^{kl}.
\end{eqnarray*}
The last inequality yields $\left\|z(t)\right\| \le \alpha_3 \bar m.$ Because $\alpha_3<1$ we obtain a contradiction. Therefore, $\overline w(t)=\widetilde w(t)$ for $t\in[0,\theta_{p+1}].$ Utilizing induction one can easily prove the uniqueness for all $t \ge 0.$ $\square$

\begin{remark} In the proof of Theorem \ref{thm2}, we make use of the contraction mapping principle to prove the exponential stability. In the literature, Lyapunov-Krasovskii functionals, LMI technique, free weighting matrix method and differential inequality technique were used to investigate the exponential stability in neural networks \cite{10,Liao02,He06}. They may also be considered in the future to prove the exponential stability in networks of the form (\ref {3}).
\end{remark}

The next section is devoted to the existence as well as the exponential stability of almost periodic solutions of (\ref{3}).

\section{Almost periodic solutions}

Let us denote by $B_0(\rr)$ the set of all bounded and continuous functions defined on $\rr.$ For $g\in B_0(\rr)$ and $ \alpha \in \rr,$ 
a translation of $g$ by $\alpha$ is a function $Q_{\alpha} g (t)= g(t+\alpha),$ $t \in \rr.$ A number $\alpha \in \rr$ is called an $\epsilon-$translation number of a function $g\in B_0(\rr)$ if $\,|| Q_{\alpha} g (t)- g(t)||< \epsilon\,$ for every $t \in \rr.$ Besides, a set $S \subset \rr$ is said to be relatively dense if there  exists a number $h > 0$ such that $[\vartheta, \vartheta+ h] \cap S \not = \emptyset $ for all $\vartheta \in \rr.$ A function $g \in B_0 (\rr)$ is said to be almost periodic, if for every positive number $\epsilon,$ there exists a relatively dense set of $\epsilon-$translation numbers of $g$ \cite{cord}.

On the other hand, an integer $k_0$ is called an $\epsilon-$almost period of a sequence $\left\{a_p\right\},$ $p\in\mathbb Z,$ of real numbers if $\left|a_{p+k_0}-a_p\right|<\epsilon$ for any $p\in\mathbb Z$ \cite{sp}. Let $\zeta_p^q = \zeta_{p+q} -\zeta_p$ and $\theta_p^q = \theta_{p+q} - \theta_p$ for all $p$ and $q.$ We call the family of sequences $\left\{\zeta_p^q\right\},$ $q \in \zz,$ equipotentially almost periodic \cite{akhmet,sp,hala} if for an arbitrary positive number $\epsilon$ there exists a relatively dense set of $\epsilon-$almost periods, common for all sequences $\left\{\zeta_p^q\right\},$ $q \in \zz.$

The following conditions are required.
\begin{itemize}
\item[\bf (C8)] The sequences $\left\{\zeta_p^q\right\},$ $q \in \zz,$ as well as the sequences $\left\{\theta_p^q\right\},$ $q \in \zz,$ are equipotentially almost  periodic;
\item [\bf (C9)] There exist positive numbers $\underline \theta$ and $\underline \zeta$ such that $\theta_{p+1} - \theta_{p} \ge \underline \theta$ and $\zeta_{p+1} - \zeta_{p} \ge \underline \zeta$ for all $p \in \mathbb Z.$ 
\end{itemize}
 
It follows from condition $(C8)$ that there exists a positive number $\bar \theta$ such that condition $(C3)$ is valid, and $|\theta_p|,|\zeta_p| \to \infty$ as $|p| \to \infty$ \cite{akhmet,sp,hala}.

The next assertion can be proved by the method of common almost periods developed in \cite{Wexler66} (see also \cite{sp,hala}). 

\begin{lem} \label{lem1} \cite{hala} Assume that $L(t)=\left\{L_{ij}(t)\right\},$ $i = 1,2,\ldots,m,$ $j=1,2,\ldots,n$ is almost periodic and the conditions $(C8),$ $(C9)$ are valid. Then, for arbitrary $\eta >0,$ $0 <\nu <\eta,$ there exist relatively dense sets of real numbers $\Omega$ and integers $Q$ such that 
\begin{enumerate}
\item[(i)] $\|L(t+\alpha) - L(t)\|< \eta,$ $t \in \rr;$
\item[(ii)] $|\zeta_p^q - \alpha| < \nu,$ $p \in \zz;$
\item[(iii)] $|\theta_p^q - \alpha| < \nu,$ $p \in \zz,$ $\alpha \in \Omega,$ $q \in Q.$	
\end{enumerate}
\end{lem}

The existence and exponential stability of the almost periodic solution of the network (\ref{3}) is mentioned in the following theorem.

\begin{thm} \label{thm3} Assume that the conditions $(C1),(C2),$ $(C4)-(C6),$ $(C8)$ and $(C9)$ are fulfilled. Then, the $SICNN$ (\ref{3}) admits a unique almost periodic solution.  If, additionally, the condition $(C7)$ is valid, then the solution is exponentially stable with exponential convergence rate $\gamma_0/2.$  
\end{thm}
{\bf Proof.} It follows from Theorem \ref{thm2} that (\ref{3}) admits a unique bounded on $\rr$ solution $u(t)=\left\{u_{ij}(t)\right\},$ $i=1,2,\ldots,m,$ $j=1,2,\ldots,m,$ which is exponentially stable provided that the condition $(C7)$ is valid. We will show that it is an almost periodic function. 

Consider the operator $\Pi$ defined by equation (\ref{theorem_operator_defn}) again. It is sufficient to verify that $\Pi u(t)$ is almost periodic, if $u(t)$ is.

Let us denote
$
\displaystyle \beta=\max_{(i,j)} \left[\frac{1}{a_{ij}} + (M+3LH) \frac{\sum_{C_{kl}\in N_r(i,j)}C_{ij}^{kl}}{a_{ij}} + \frac{4LH^2\sum_{C_{kl}\in N_r(i,j)}C_{ij}^{kl}}{1-e^{-a_{ij} \underline{\theta}}}  \right].
$
Fix an arbitrary positive number $\epsilon.$ Because $\Pi u$ is uniformly continuous, there exists a positive number $\eta$ satisfying $\displaystyle \eta <   \frac{\underline{\theta}}{5}$ and $\displaystyle \eta \leq \frac{\epsilon}{3\beta}$ such that if $\left|t'-t''\right|<4\eta,$ then 
\begin{eqnarray} \label{SICNN_ap_proof_3}
\left\|\Pi u(t') - \Pi u (t'')\right\| < \displaystyle \frac{\epsilon}{3}.
\end{eqnarray} 
Next, we take into account a number $\nu$ with $0<\nu<\eta$ such that $\left\|u(t')-u(t'')\right\|<\eta$ whenever $\left|t'-t''\right|<\nu,$ and let $\alpha$ and $q$ be numbers as mentioned in Lemma \ref{lem1} such that $\alpha$ is an $\eta-$translation number for $u(t).$

Assume that $t\in (\theta_{\overline{p}} + \eta, \theta_{\overline{p}+1} - \eta)$ for some $\overline p \in \mathbb Z.$ Making use of the equation
\begin{eqnarray*}
&& \displaystyle(\Pi u(t+\alpha))_{ij} - (\Pi u(t))_{ij} =  -\int_{-\infty}^{t} e^{-a_{ij}(t-s)}\sum_{C_{kl} \in N_r(i,j)} C_{ij}^{kl} \\
&& \times \displaystyle \Big[ f(u_{kl(s +\alpha)},u_{kl\gamma(s+\alpha)})u_{ij}(s +\alpha)   - f(u_{kls},u_{kl\gamma(s)})u_{ij}(s)  \Big]ds \\
&& + \displaystyle \int_{-\infty}^{t} e^{-a_{ij}(t-s)} \left[ L_{ij}(s+\alpha)- L_{ij}(s)   \right] ds,
\end{eqnarray*}
we obtain that
\begin{eqnarray} \label{SICNN_ap_proof_1}
\begin{array}{l} 
\displaystyle \left|(\Pi u(t+\alpha))_{ij} - (\Pi u(t))_{ij}\right| \le \int_{-\infty}^{t} e^{-a_{ij}(t-s)}\sum_{C_{kl} \in N_r(i,j)} C_{ij}^{kl} M \left|u_{ij}(s+\alpha)-u_{ij}(s)\right|ds \\
 + \displaystyle\int_{-\infty}^{t} e^{-a_{ij}(t-s)} \sum_{C_{kl} \in N_r(i,j)} C_{ij}^{kl} LH \left\| u_{kl(s+\alpha)}-u_{kls} \right\|_0 ds \\
 + \displaystyle\int_{-\infty}^{t} e^{-a_{ij}(t-s)} \sum_{C_{kl} \in N_r(i,j)} C_{ij}^{kl} LH \left\| u_{kl\gamma(s+\alpha)}-u_{kl\gamma(s)} \right\|_0 ds \\
 + \displaystyle\int_{-\infty}^{t} e^{-a_{ij}(t-s)}  \left|L_{ij}(s+\alpha)- L_{ij}(s)\right| ds.
\end{array}
\end{eqnarray}
According to Lemma \ref{lem1}, $(i),$ the inequality 
\begin{eqnarray*}
\int_{-\infty}^{t} e^{-a_{ij}(t-s)}  \left|L_{ij}(s+\alpha)- L_{ij}(s)\right| ds < \frac{\eta}{a_{ij}}
\end{eqnarray*}
is valid. Moreover, since $\alpha$ is an $\eta-$translation number for $u(t),$ one can confirm that
\begin{eqnarray*}
\int_{-\infty}^{t} e^{-a_{ij}(t-s)} \sum_{C_{kl} \in N_r(i,j)} C_{ij}^{kl} M \left|u_{ij}(s+\alpha)-u_{ij}(s)\right|ds <  \frac{M \eta}{a_{ij}} \sum_{C_{kl} \in N_r(i,j)} C_{ij}^{kl}
\end{eqnarray*}
and
\begin{eqnarray*}
\displaystyle\int_{-\infty}^{t} e^{-a_{ij}(t-s)} \sum_{C_{kl} \in N_r(i,j)} C_{ij}^{kl} LH \left\| u_{kl(s+\alpha)}-u_{kls} \right\|_0 ds <  \frac{LH \eta}{a_{ij}}\sum_{C_{kl} \in N_r(i,j)} C_{ij}^{kl}.
\end{eqnarray*}
On the other hand, we have
\begin{eqnarray}\label{SICNN_ap_ineq2}
\begin{array}{l}
 \displaystyle \int_{-\infty}^{t} e^{-a_{ij}(t-s)} \sum_{C_{kl} \in N_r(i,j)} C_{ij}^{kl} LH \left\| u_{kl\gamma(s+\alpha)}-u_{kl\gamma(s)} \right\|_0 ds \\
 \le \displaystyle\int_{\theta_{\overline{p}}+\eta}^{t} e^{-a_{ij}(t-s)} \sum_{C_{kl} \in N_r(i,j)} C_{ij}^{kl} LH \left\| u_{kl\gamma(s+\alpha)}-u_{kl(\gamma(s)+\alpha)} \right\|_0 ds \\
 + \displaystyle \sum_{\lambda=0}^{\infty} \int_{\theta_{\overline{p}-\lambda-1}+\eta}^{\theta_{\overline{p}-\lambda}-\eta} e^{-a_{ij}(t-s)} \sum_{C_{kl} \in N_r(i,j)} C_{ij}^{kl} LH \left\| u_{kl\gamma(s+\alpha)}-u_{kl(\gamma(s)+\alpha)} \right\|_0 ds \\
 + \displaystyle \sum_{\lambda=0}^{\infty} \int_{\theta_{\overline{p}-\lambda}-\eta}^{\theta_{\overline{p}-\lambda}+\eta} e^{-a_{ij}(t-s)} \sum_{C_{kl} \in N_r(i,j)} C_{ij}^{kl} LH \left\| u_{kl\gamma(s+\alpha)}-u_{kl(\gamma(s)+\alpha)} \right\|_0 ds \\
 + \displaystyle \int_{-\infty}^{t} e^{-a_{ij}(t-s)} \sum_{C_{kl} \in N_r(i,j)} C_{ij}^{kl} LH \left\| u_{kl(\gamma(s)+\alpha)}-u_{kl\gamma(s)} \right\|_0 ds.
 \end{array}
\end{eqnarray}
For any $p\in\mathbb Z,$ if $s\in (\theta_p+\eta, \theta_{p+1}-\eta),$ then one can show by using Lemma $\ref{lem1}, (iii)$ that the number $s +\alpha$ belongs to the interval $(\theta_{p + q},\theta_{p + q+1})$ so that  
$$
\left\| u_{\gamma(s+\alpha)} -u_{\gamma(s)+\alpha} \right\|_0 = \max_{\kappa\in [-\tau,0]} \left| u(\kappa+\zeta_{p+q}) - u(\kappa+\zeta_p+\alpha) \right| < \eta,
$$
since $\left| (\kappa+\zeta_{p+q})-(\kappa+\zeta_p+\alpha)  \right| = \left| \zeta^q_p - \alpha \right|<\nu$ by Lemma \ref{lem1}, $(ii).$ 
Besides, the inequality 
$$
\displaystyle \sum_{\lambda=0}^{\infty} \int_{\theta_{\overline{p}-\lambda}-\eta}^{\theta_{\overline{p}-\lambda}+\eta} e^{-a_{ij}(t-s)} ds \le 2\eta \displaystyle \sum_{\lambda=0}^{\infty} e^{-a_{ij}\underline{\theta}\lambda}=\frac{2\eta}{1-e^{-a_{ij}\underline{\theta}}}
$$ 
is valid.
Therefore, (\ref{SICNN_ap_ineq2}) yields
\begin{eqnarray*}
&& \displaystyle \int_{-\infty}^{t} e^{-a_{ij}(t-s)} \sum_{C_{kl} \in N_r(i,j)} C_{ij}^{kl} LH \left\| u_{kl\gamma(s+\alpha)}-u_{kl\gamma(s)} \right\|_0 ds \\
&& <  LH \eta \left( \frac{2\sum_{C_{kl} \in N_r(i,j)} C_{ij}^{kl}}{a_{ij}} + \frac{4H\sum_{C_{kl} \in N_r(i,j)} C_{ij}^{kl}}{1-e^{-a_{ij}\underline{\theta}}}  \right).
\end{eqnarray*}
It can be verified by means of (\ref{SICNN_ap_proof_1}) that 
$$
 \displaystyle \left|(\Pi u(t+\alpha))_{ij} - (\Pi u(t))_{ij} \right| 
  < \eta\left[\frac{1}{a_{ij}} + (M+3LH) \frac{\sum_{C_{kl}\in N_r(i,j)}C_{ij}^{kl}}{a_{ij}} + \frac{4LH^2\sum_{C_{kl}\in N_r(i,j)}C_{ij}^{kl}}{1-e^{-a_{ij} \underline{\theta}}} \right]. 
$$
Hence, 
\begin{eqnarray} \label{SICNN_ap_proof_4}
\left\|\Pi u (t+\alpha) - \Pi u(t)\right\|<\beta \eta \le \epsilon/3
\end{eqnarray} 
for each $t$ that belongs to the intervals $(\theta_{p} + \eta, \theta_{p+1} - \eta),$ $p\in \mathbb Z.$

The inequality $\eta < \underline{\theta}/5$ ensures that $t+3\eta \in (\theta_{p} + \eta, \theta_{p+1} - \eta)$ if $\left| t-\theta_p  \right| \le \eta.$
Now, by means of the inequalities (\ref{SICNN_ap_proof_3}) and (\ref{SICNN_ap_proof_4}) we attain for $\left| t-\theta_p \right| \le \eta,$ $p\in\mathbb Z,$ that
\begin{eqnarray*}
&& \left\|\Pi u(t+\alpha) - \Pi u(t)\right\| \le \left\|\Pi u(t+\alpha) - \Pi u(t+\alpha+3\eta)\right\| \\
&&+ \left\|\Pi u(t+\alpha+3\eta) - \Pi u(t+3\eta)\right\| + \left\|\Pi u(t+3\eta) - \Pi u(t)\right\| \\
&&<\epsilon.
\end{eqnarray*}
The last inequality implies that $\alpha$ is an $\epsilon-$translation number of $\Pi u(t).$ Consequently, the $SICNN$ (\ref{3}) admits a unique almost periodic solution. 
 $\square$

\begin{remark}
The Bohr definition of almost periodicity is also suitable for the application of Lyapunov functional method and the technique of Young inequality \cite{Jiang05} to show the existence, uniqueness and exponential stability of almost periodic solutions in $CNNs.$
\end{remark}

\section{An example}

Consider the sequence $\theta=\left\{\theta_p\right\}$ defined as $\displaystyle \theta_p = p + \frac{1}{4}|\sin(p) - \cos(p \sqrt 2)|,$ $p\in\mathbb Z.$ Utilizing the technique provided in \cite{sp,akhmet54}, one can verify that the sequences  $\left\{\theta_p^q\right\},$ $q \in \mathbb Z,$ are equipotentially  almost periodic. We take the function $\gamma(t)$ with $\zeta_p=\theta_p.$ One can confirm that the conditions  $(C3)$ and $(C9)$ hold with $\bar{\theta}=3/2$ and $\underline{\theta}=1/2,$ respectively.
 
Let us take into account the $SICNN$
\begin{eqnarray}\label{pca_example}
\frac{dx_{ij}}{dt} = - a_{ij} x_{ij} - \sum_{C_{kl} \in N_1(i,j)} C_{ij}^{kl}f(x_{kl}(\gamma(t)-\tau))x_{ij} +  L_{ij}(t),
\end{eqnarray} 
in which $i,j = 1,2,3,$ $f(s)=\displaystyle \frac{s^2}{2}$ if $|s| \le 0.1,$ $f(s)=0.005$ if $|s|>0.1,$ $\tau=0.3,$
$$\left( \begin{array}{ccc}
a_{11}&a_{12}&a_{13} \\
a_{21}&a_{22}&a_{23} \\
a_{31}&a_{32}&a_{33} \end{array} \right)= \left( \begin{array}{ccc}
9&3&5 \\
6&5&4 \\
3&12&9 \end{array} \right),$$  
$$\left( \begin{array}{ccc}
C_{11}&C_{12}&C_{13} \\
C_{21}&C_{22}&C_{23} \\
C_{31}&C_{32}&C_{33} \end{array} \right)= \left( \begin{array}{ccc}
0.08&0.01&0.02 \\
0.05&0.03&0.06 \\
0.04&0.07&0.02 \end{array} \right),$$
\begin{eqnarray*}
&& \left( \begin{array}{ccc}
L_{11}(t)&L_{12}(t)&L_{13}(t) \\
L_{21}(t)&L_{22}(t)&L_{23}(t) \\
L_{31}(t)&L_{32}(t)&L_{33}(t) \end{array} \right) \\
&& = \left( \begin{array}{ccc}
0.1\cos(t) + 0.2\sin(\sqrt{2}t)&0.2\cos(\pi t) + 0.1\sin(\sqrt{2}t)&0.15\cos(2t)- 0.12\cos(\pi t)\\
0.15\cos(3t) - 0.1\sin(\pi t)&0.2\cos(t) - 0.15\sin(\sqrt{2}t)&0.1\sin(t)+ 0.2\cos(\sqrt{3}t)\\
0.2\cos(\sqrt{2}t) + 0.14\sin(\pi t)&0.2\cos(\sqrt{2}t) + 0.1\sin(t)&0.15\cos(\sqrt{2}t)- 0.13\cos(4t) \end{array} \right).
\end{eqnarray*} 

One can calculate that 
$ \sum_{C_{kl} \in N_1(1,1)} C_{11}^{kl} =0.17,$ $\sum_{C_{kl} \in N_1(1,2)} C_{12}^{kl} =  0.25,$ $\sum_{C_{kl} \in N_1(1,3)} C_{13}^{kl} = 0.12,$
$\sum_{C_{kl} \in N_1(2,1)} C_{21}^{kl} = 0.28,$ $\sum_{C_{kl} \in N_1(2,2)} C_{22}^{kl} = 0.38,$ $\sum_{C_{kl} \in N_1(2,3)} C_{23}^{kl} = 0.21,$ $\sum_{C_{kl} \in N_1(3,1)} C_{31}^{kl} = 0.19,$ $\sum_{C_{kl} \in N_1(3,2)} C_{32}^{kl} = 0.27,$ $\sum_{C_{kl} \in N_1(3,3)} C_{33}^{kl} = 0.18.$ The conditions $(C5)-(C7)$ are valid for (\ref{pca_example}) with $\gamma_0=3,$ $\mu=0.38,$ $\bar{c}=\bar{d}=0.25/3,$ $M=0.005,$ $L=0.1,$ $\bar{L}=0.35,$ $\bar{l}=0.34/3.$ According to Theorem \ref{thm3}, the network (\ref{pca_example}) has a unique almost periodic solution, which is exponentially stable with the rate of convergence $3/2.$

\begin{figure}[ht] 
\centering
\includegraphics[width=15.2cm]{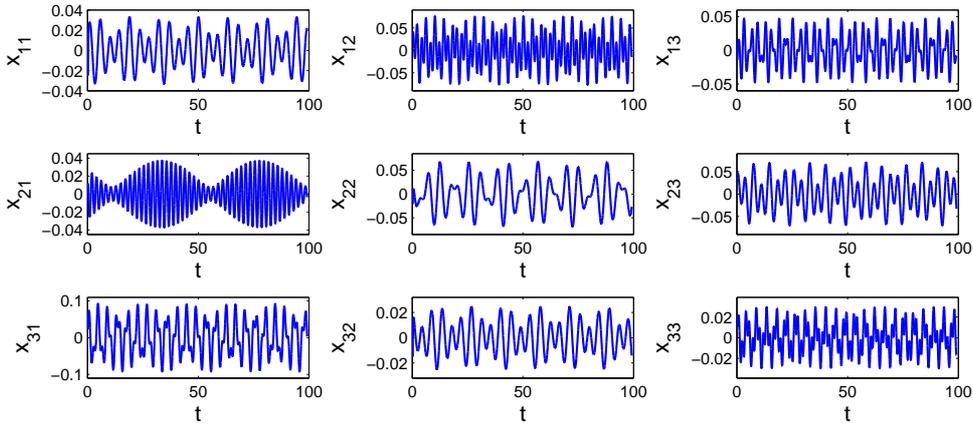}
\caption{The unique almost periodic solution of $SICNN$ (\ref{pca_example}).}
\label{pca_fig1}
\end{figure}

Consider the constant function $\phi(t)=\left\{\phi_{ij}(t)\right\}$ such that $\phi_{11}(t)=-0.025,$ $\phi_{12}(t)=0.036,$ $\phi_{13}(t)=-0.014,$ $\phi_{21}(t)=0.012,$ $\phi_{22}(t)=-0.021,$ $\phi_{23}(t)=0.042,$ $\phi_{31}(t)=0.023,$ $\phi_{32}(t)=-0.015,$ $\phi_{33}(t)=0.012.$ We depict in Figure \ref{pca_fig1} the solution $x(t)=\left\{x_{ij}(t)\right\}$ of (\ref{pca_example}) with $x(t)=\phi(t),$ $t\le \sigma= \theta_0= \displaystyle \frac{1}{4}.$
 Figure \ref{pca_fig1} supports the result of Theorem \ref{thm3} such that the represented solution converges to the unique almost periodic solution of $SICNN$ (\ref{pca_example}).

\section{Conclusion} 
 
In this paper, we investigate the existence as well as the exponential stability of almost periodic solutions in a new model of $SICNNs.$ The usage of the functional response on alternate (advanced-delayed) type of piecewise constant arguments is the main novelty of our study, and it is useful for the investigation of a large class of neural networks. An illustrative example is provided to show the effectiveness of the theoretical results.

Our approach concerning exponential stability may be used in the future to investigate synchronization of chaos in coupled neural networks and control of chaos in large communities of neural networks with piecewise constant argument.

Differential equations with functional response on piecewise constant argument can be applied for the development of other kinds of recurrent networks such as Hopfield and Cohen-Grossberg neural networks \cite{Hopfield84,Cohen93} and others. This will provide new opportunities for the analysis and applications of neural networks.

\section*{Acknowledgments}

The authors wish to express their sincere gratitude to the referees for the helpful criticism and valuable suggestions, which helped to improve the paper significantly. 

The second author is supported by the 2219 scholarship programme of T\"{U}B\.{I}TAK, the Scientific and Technological Research Council of Turkey.

\end{document}